\documentclass[12pt]{article}
\usepackage{mathrsfs}
\usepackage{amsfonts}
\usepackage{amssymb, amsmath,mathrsfs,amsfonts}

\numberwithin{equation}{section}

\begin{document}
\title {\large SPECTRA OF LINEAR FRACTIONAL COMPOSITION OPERATORS ON
$H^2(B_N)$\footnotetext{ {\it 2000 Mathematics Subject
Classification.} Primary: 47B33; Secondary 32A35.} \footnotetext{
{\it Key words and phrases.} Composition operators, Hardy space,
spectra, linear fractional maps.}}
\author{\normalsize  LIANGYING JIANG\footnote{Liangying Jiang is
supported by  Shanghai Education Research and Innovation Project
(No. 10YZ185)
 and by Shanghai University Research Special Foundation  for Outstanding Young Teachers
 (No. sjr09015)}
 \quad ZHIHUA  CHEN\footnote{ Zhihua Chen  is
supported by the National Natural Science Foundation of China (No.
10871145) and Doctoral Program Foundation of the Ministry of
Education of China (No. 20090072110053)}}
\date{}
\maketitle \leftskip=10mm \rightskip=10mm \noindent{\small ABSTRACT.
We characterize the spectra of  composition operators on the Hardy
space $H^2(B_N)$, when the symbols are  elliptic or  hyperbolic
linear fractional self-maps of $B_N$. Therefore,  combining with the
result obtained by Bayart \cite{B10}, the spectra of all linear
fractional composition operators on $H^2(B_N)$ are completely
determined.}

\leftskip=0mm \rightskip=0mm

\section*{\large  1 Introduction}

Let $B_N$ denote the unit ball of $\mathbb{C}^N$, and let $\varphi:
B_N\rightarrow B_N$ be an analytic map. In this paper, we consider
the composition operator $C_\varphi$ defined by
$C_\varphi(f)=f\circ\varphi$, acting on the classical Hardy space
$H^2(B_N)$. For $N>1$, some authors gave the examples  of unbounded
composition operators on $H^2(B_N)$ (see Section 3.5 of \cite{CM1}),
which exhibit surprisingly different behaviors with the case of one
complex variable. So many properties of composition operators in
several variables are not easily managed.

However, if the symbol $\varphi$ is a linear fractional self-map of
$B_N$, that is,  $$\varphi(z)=\frac{Az+B}{<z, C>+d}$$ for  an
$N\times N$ matrix
 $A=(a_{jk})$, two column vectors $B=(b_j)$, $C=(c_j)$  of $N\times 1$ and $d\in
 \mathbb{C}$, such that $\varphi$ maps the unit ball $B_N$ into
 itself,   applying Wogen's criterion \cite{Wo}, Cowen and MacCluer \cite{CM2} proved that  $C_\varphi$ is
 bounded on $H^2(B_N)$.
 In recent years, linear fractional maps  of
 the  ball and their  composition operators on some analytic function spaces  have been
 developed from various aspects, for example, geometric properties
 \cite{BB} and
 classification of semigroups \cite{BCD} for  linear fractional maps; cyclic
 behavior (see \cite{B}, \cite{B10}, \cite{JO08})  and  essential
 normality (see  \cite{JO}, \cite{MW}, \cite{ZY09}) of
their composition operators. Similar to the case of the unit disk
$D$, linear
 fractional self-maps of $B_N$  give rich examples to exhibit
the  complexity of  properties of the associated composition
 operators.  Moreover, motivated by the fact that the linear fractional
 model in  $D$ provides a useful tool to deal with composition
 operators induced by general maps, some authors began to study a linear fractional model for general maps of
 the unit ball, see  \cite{B08} and \cite{CC}.

In this paper, we are interested in the spectra of linear
 fractional composition operators on  $H^2(B_N)$, this is a
 continuation of \cite{B10}. As a part of operator theory, the
 structures  of spectra of  composition operators are   very important. Several authors have used this tool to investigate the
 cyclicity of composition operators (see  \cite{B10}, \cite{GM}).
 In the unit disk, the spectra of invertible and compact composition
 operators on $H^2(D)$ have been described completely, while understanding the spectra of general composition
 operators is not easy. When  $\varphi$ is a linear fractional self-map  of
 $D$, the spectrum of $C_\varphi$ on $H^2(D)$ (see \cite{Co}) is   well
 known. Recently,  Higdon \cite{Hig} completely gave the spectrum of
$C_\varphi$ on the Dirichlet space $\mathcal {D}$.
 Moreover, Hurst \cite{Hu} obtained some  results for the spectrum of  $C_\varphi$ on some weighted Hardy
 spaces. These results not only show the diversity of spectra of composition operators, and much of the spectral
 information depends on the behavior of $\varphi$ near the Denjoy-Wollf point.  From the table of  \cite{Co}, one has the following
 results.
\\ \\ {\bf Theorem A.} \begin{em}  Suppose that $\varphi$ is a linear fractional
self-map of $D$, not
 automorphism.
The spectrum of $C_\varphi$ on $H^2(D)$ can be described as
 follows:

(i) If $\varphi$  has an interior fixed point $a$, there exists two
cases. When $\varphi$ fixes a boundary point  $\tau$, then
$$\sigma(C_\varphi)=\{\lambda: |\lambda|\le
\varphi'(\tau)^{-1/2}\}\cup\{1\}.$$ Otherwise, $C_\varphi^n$ is
compact for some positive integer  $n$ and
$$\sigma(C_\varphi)=\{\varphi'(a)^k: k=0, 1, \ldots\}\cup\{0\}.$$

(ii) If $\varphi$ is parabolic, by the Cayley transformation,
$\varphi$ is conjugated to a map $\psi(z)=z+t$ on the upper-half
plane, then
$$\sigma(C_\varphi)=\{e^{\beta t}:\beta\le 0\}\cup\{0\}.$$

(iii) If $\varphi$ is hyperbolic with the  Denjoy-Wollf point
$\tau\in \partial D$, then $$\sigma(C_\varphi)=\{\lambda:
|\lambda|\le \varphi'(\tau)^{-1/2}\}.$$ \end{em} \\ \par In the unit
ball, only the  spectra of automorphism-induced composition
operators and compact operators have been  determined (see Chapter 7
in \cite{CM1}).  For  other cases, less results have been obtained,
even for the symbols being  linear fractional self-maps of $B_N$.
Recently, Jury \cite{J2} has computed the spectral radii of linear
fractional composition operators on $H^2(B_N)$, soon later, the
spectra of composition operators were
 characterized by Bayart \cite{B10}  when their  symbols are
parabolic linear fractional maps of $B_N$. In order to illustrate
these results, we recall some definitions for the  Denjoy-Wollf
points and the classification of linear fractional maps of $B_N$. In
this paper, we denote by LFM($B_N$) the set of linear fractional
maps of $B_N$.

The Denjoy-Wollf Theorem of the unit  ball is the following (see
\cite{M83}).
\\ \\ {\bf Theorem B.} \begin{em} Let $\varphi\in {\mbox LFM}(B_N)$ with no fixed
points in $B_N$.
 There exists a unique point $\tau\in\partial B_N$ such that $\varphi(\tau)=\tau$  and $<d
\varphi_{\tau}(\tau), \tau>=\alpha$ with $0<\alpha\le 1$.\end{em} \\
\\ The point $\tau$ is called the Denjoy-Wollf point of $\varphi$ and
$\alpha$ is the boundary dilation coefficient. In particular, we
call such a map hyperbolic when $\alpha\in (0,1)$ and parabolic when
$\alpha=1$. Otherwise, if $\varphi$ has a fixed point in $B_N$, then
we call it elliptic.

Now, we have the  following  result for   spectral radii of linear
fractional composition operators on $H^2(B_N)$ (see \cite{J2}).
\\ \\ {\bf Theorem C.} \begin{em} Let $\varphi\in {\mbox LFM}(B_N)$. The spectral
radius of $C_\varphi$ acting on  $H^2(B_N)$ is $1$ if $\varphi$ is
elliptic; if $\varphi$ is non-elliptic with the boundary dilation
coefficient $\alpha$, then the spectral radius is
$\alpha^{-N/2}$.\end{em} \\ \par For a parabolic linear fractional
self-map $\varphi$ of $B_N$, considering  its conjugating map $\psi$
on the Siegel half-plane $H_N$, Bayart \cite{B10} obtained a normal
form of  $\psi$ and gave a classification for this normal form.
Based on this classification,
 he completely  computed  the spectrum of $C_\varphi$ on $H^2(B_N)$.

In this paper, we continue  to investigate the spectra of linear
fractional composition operators on  $H^2(B_N)$ induced by the
remaining maps, and  we organize it as follows. In Section 2, we
first look for a classification for elliptic linear fractional
self-maps of $B_N$, the idea comes from the classification of
elliptic semigroups (see \cite{BCD}) and the geometric
classification according to their  fixed point sets (see \cite{BB}).
Different to the case of the disk, elliptic linear fractional  maps
of $B_N$ contain three cases. Therefore, we will prove that the
spectra of the  associated composition operators  have three
different structures. For a general map $\varphi$, univalent, not
automorphism, which fixes a point in $D$ or $B_N$, the spectrum of
$C_\varphi$ on many function spaces has been studied (see
\cite{CM94}, \cite{LP},  \cite{MS}, \cite{YZ}, \cite{Zheng}). It is
easy to see that  the proofs of these results follow the same
pattern. We shall use an analogue approach to deal with the spectrum
of $C_\varphi$ on  $H^2(B_N)$, when the symbol  $\varphi$ is an
elliptic linear fractional map with  only one boundary fixed  point.

In Section 3, we devote to  hyperbolic linear fractional self-maps
$\varphi$ of $B_N$. It is well  known  that  $\varphi$ has only one
or two fixed points on the boundary $\partial B_N$ (see \cite{BB}).
First, we will obtain a simple form for each conjugation map $\psi$
of $\varphi$ on the Siegel-half plane  $H_N$. Applying these forms,
we can completely determine the spectra of hyperbolic composition
 operators on $H^2(B_N)$. In particular, when $\varphi$ fixes only one boundary point, the spectrum of $C_\varphi$
 has a similar structure with the corresponding case in the disk,
 which has been proved in the first author's doctoral thesis
 \cite{J08}. For the case $\varphi$  fixing two boundary  points, Xu
 and Deng \cite{XD} recently used the method of Bayart  \cite{B10} to  obtain the spectrum of $C_\varphi$ on
 $H^2(B_N)$. Thus,  the    spectra of all linear fractional composition
 operators on $H^2(B_N)$ are completely characterized.

 \section*{\large 2 Spectra of elliptic composition operators}

 In this section, we deal with elliptic linear fractional self-maps
of $B_N$ and the spectra of the associated composition operators.
First, we need a classification for these maps.

Recall that a slice $S$ is a non-empty subset of $B_N$ of the form
$S=B_N\cap V$, where $V$ is a one-dimensional affine subspace of
$\mathbb{C}^N$. A $p$-dimensional slice of $B_N$ is the non-empty
intersection between $B_N$ and a $p$-dimensional affine subspace of
$\mathbb{C}^N$ with $p\ge 0$. By Herv\'{e}'s theorem (see \cite{Ab}
or \cite{Ru}), if $\varphi\in {\mbox LFM}(B_N)$ has a non-empty
fixed points set in $B_N$ then such a set is a $p$-dimensional slice
of $B_N$.

In \cite{BB}, Bisi and Bracci gave a geometric classification for
linear fractional self-maps of $B_N$ based on their fixed point
sets.  In this classification, elliptic linear fractional maps
contain three cases. For elliptic semigroup of linear fractional
self-maps of $B_N$, Bracci et al \cite{BCD} also classified this
semigroup into three different cases. Surprisedly, we find that
there is a close connection between two classifications and  obtain
the following result. Let $$L_U(\varphi,
z_0)=\bigoplus\limits_{|\alpha|=1} \mbox{ker}(d\varphi_{z_0}-\alpha
I)^N$$ for $\varphi\in \mbox{Hol}(B_N, B_N)$ with a fixed point
$z_0\in B_N$, which is called the unitary space of $\varphi$ at the
point $z_0$, the dimension of it is called the unitary index of
$\varphi$ (see \cite{BCD}).
\\ \\ {\bf Theorem 2.1.}\begin{em} Let $\varphi\in {\mbox LFM}(B_N)$ be elliptic
with a fixed point $z_0\in B_N$  and let $p=\mbox{dim}\ L_U(\varphi,
z_0)$.

(1) If $p>0$, then $\varphi$ is conjugated to a map $\psi$ with
$$\psi(z',z'')=(Uz', Az''),\qquad (z',z'')\in \mathbb{C}^p\times\mathbb{C}^{N-p}\cap B_N,$$
where $U$ is a unitary diagonal matrix of  $\mathbb{C}^{p\times p}$
and $A$ is a matrix of order $N-p$ with $||A||< 1$.

(2) If $p=0$, $\varphi$ will fix at most one boundary point. In this
case:

\quad (i) When $\varphi$ has no  boundary fixed point, then
$\varphi$ is conjugated to  a map of the form $$\psi(w)=Aw$$ defined
on a complex ellipsoid
$$\triangle_1=\{w=(w_1, w')\in \mathbb{C}\times\mathbb{C}^{N-1}:
\frac{1}{r^2}|w_1-\sqrt{r^2-1}|^2+|w'|^2<r^2\}$$ for some $r\ge 1$,
where $A$ is a matrix of $\mathbb{C}^{N\times N}$ with $||A||<1$.

\quad (ii) When $\varphi$ has only one boundary fixed point, then
$\varphi$ is conjugated to  $$\psi(w)=Aw$$  on a half-plane
$$\triangle_2=\{w=(w_1, w')\in \mathbb{C}\times\mathbb{C}^{N-1}:
\mbox{Re}\, 2w_1>|w'|^2-1\},$$ where $A$ is a matrix of order $N$
 with $||A||<1$.\end{em}
\\ \par {\bf Proof.} We will use a similar argument as that of Theorem 3.2
and Corollary 3.3 in \cite{BCD} to obtain these results. Especially,
the proof of (1) can be seen in \cite{BCD} or Proposition 3.8 of
\cite{J08}. Here, we only give a proof for (2).

If $p=0$, it is clear that $\varphi$ fixes only one interior point
$z_0$. Up to conjugation with automorphisms of $B_N$, we may assume
that $z_0=0$. We claim that $\varphi$ fixes at most one boundary
point. Otherwise, if $\varphi$ fixes two points on $\partial B_N$,
by Theorem 3.2 of \cite{BB}, $\varphi$ is conjugated to a map which
has a hyperbolic automorphism (not identity) at first coordinate.
Thus, $\varphi$ has no fixed point in $B_N$, which contradicts with
$\varphi(0)=0$.  When $\varphi$ has more than two boundary fixed
points, similar to the argument in the proof of Theorem 3.1 in
\cite{BB}, $\varphi$ will fix a $q$-dimensional slice of $B_N$ with
$q\ge 1$. This means that the unitary index $p$ of $\varphi$ must be
larger than $q$, contradicting to the hypothesis $p=0$. So that the
claim is true.

Now,  $\varphi(0)=0$ gives  $$\varphi(z)=\frac{Az}{<z, C>+1}$$ for
some $A\in \mathbb{C}^{N\times N}$  and $C\in \mathbb{C}^N$. Since
$\varphi$ maps $B_N$ into $B_N$, we have $|C|<1$. A computation
shows that $d\varphi_0=A$. Thus  $p=0$ implies $||A||<1$, that is,
$A^\ast-I$ is invertible. Immediately, there exists a vector $V\in
\mathbb{C}^N$  such that $(A^\ast-I)V=C$  and $\delta:=|V|\le 1$.
Conjugating $\varphi$ by a unitary map $U$ with $U^\ast V=\delta
e_1$, where $e_1=(1,0,\ldots,0)\in\partial B_N$, we obtain
$$\widetilde{\varphi}(z)=\frac{A_1z}{\delta<z, (A_1^\ast-I)e_1>+1}$$
with $A_1=U^\ast A U$ and $||A_1||<1$. Define
$$\sigma(z)=\frac{z}{-\delta z_1+1},\qquad z=(z_1, z')\in \mathbb{C}\times\mathbb{C}^{N-1}\cap
B_N,$$  it is an one-to-one holomorphic linear fractional map from
$B_N$ onto $\Omega:=\sigma(B_N)$. Consequently, we deduce that
$$\sigma\circ \widetilde{\varphi}(z)=A_1 \sigma(z), \qquad z\in B_N,$$
and $$\Omega=\{w=(w_1, w')\in\mathbb{C}\times\mathbb{C}^{N-1}:
|1+\delta w_1|^2>|w'|^2 \}.$$

If $\varphi$ fixes no boundary  point, then $\delta<1$. Otherwise,
$\delta=1$ implies that
$$\widetilde{\varphi}(z)=\frac{A_1z}{<z, (A_1^\ast-I)e_1>+1}$$
has a boundary  fixed  point $e_1$.  So $\varphi$ also fixes a
boundary point, which is a contradiction. If we  set
$r=(1-\delta^2)^{-1/2}$, then $r\ge 1$ and
$$\Omega=\triangle_1=\{(w_1, w')\in \mathbb{C}\times\mathbb{C}^{N-1}:
\frac{1}{r^2}|w_1-\sqrt{r^2-1}|^2+|w'|^2<r^2\}.$$

If $\varphi$ fixes only one boundary point, then
$\widetilde{\varphi}$ also has only one  fixed point $\tau\in
\partial B_N$. It follows that $$A_1 \sigma(\tau)=\sigma\circ
\widetilde{\varphi}(\tau)=\sigma(\tau),$$ that is, $\sigma(\tau)$ is
a fixed point of the matrix $A_1$. since $||A_1||<1$, we see that
$\sigma(\tau)=0$ or $\sigma(\tau)=\infty$. However, $\sigma$ is an
one-to-one map from $B_N$ onto $\Omega$ and $\sigma(0)=0$. This
implies $\sigma(\tau)=\infty$. Thus, we have $-\delta\tau_1+1=0$ and
$\delta=\frac{1}{\tau_1}\ge 1$, where $\tau_1$ denotes the first
coordinate of $\tau$. Which combining with $\delta\le 1$ yields
$\delta=1$. Hence, a calculation shows that
$$\Omega=\triangle_2=\{(w_1, w')\in \mathbb{C}\times\mathbb{C}^{N-1}:
\mbox{Re}\, 2w_1>|w'|^2-1\},$$ which is a domain similar to the
Siegel half-plane. \ \ $\Box$

Next, we will apply this classification to compute the spectrum of
$C_\varphi$  on $H^2(B_N)$, when the symbol $\varphi$ is  a elliptic
linear fractional self-map of $B_N$. In the first case of Theorem
2.1, $\varphi$ is conjugated to a  map $\psi(z',z'')=(Uz', Az'')$.
Let $\Lambda_\varphi=\{\lambda_1,\ldots,\lambda_N\}$ be the union of
the spectrum of $U$ and $A$, and  write $\Lambda_\varphi^\alpha$ for
$\lambda_1^{\alpha_1}\cdots \lambda_N^{\alpha_N}$, where $\alpha$ is
a  multi-index of $\mathbb{N}^N$. In this case, we have the
following result for the spectrum of $C_\varphi$, the proof is
inspired by that of \cite{B10} in computing the spectra of parabolic
composition operators on $H^2(B_N)$.
\\ \\ {\bf Theorem 2.2.}\begin{em}  Let $\varphi\in {\mbox LFM}(B_N)$ be
elliptic, non-automorphism, and let $p=\mbox{dim}\ L_U(\varphi,
z_0)$ for its fixed point $z_0\in B_N$. When $p>0$, the spectrum of
$C_\varphi$ on $H^2(B_N)$ is

(i) a union of circles  if  $U$ has an irrational unimodular
eigenvalue: $$\sigma(C_\varphi)=\bigcup\limits_{\alpha\in
\mathbb{N}^N}\Lambda_\varphi^\alpha \mathbb{T}\cup\{0\};$$

(ii) $\sigma(C_\varphi)=\bigcup\limits_{\alpha\in
\mathbb{N}^N}\Lambda_\varphi^\alpha\cup\{0\}$ if all  eigenvalues of
$U$ are rational.\end{em}
\\ \par {\bf Proof.}  By   Theorem 2.1, we see that $\varphi$ is conjugated
to a map  with the form $$\psi(z,w)=(Uz, Aw), \qquad (z, w)\in
\mathbb{C}^p\times\mathbb{C}^{N-p}\cap B_N,$$ where
$U=\mbox{diag}(e^{i\theta_1}, \ldots, e^{i\theta_p})$ and $||A||<1$.
Up to conjugation by a unitary map, we may assume that $A$ is
upper-triangular with the diagonal entries $\lambda_1,\ldots,
\lambda_{N-p}$, which  obviously are the eigenvalues of $A$. Next,
we will determine the spectrum of $C_\psi$.

For each $i\in\{1,\ldots, N-p\}$, let $w(i)=(w_{i, 1},\ldots, w_{i,
N-p})^\top$ be a non-zero eigenvector of $A^\top$ associated to
$\lambda_i$, i.e. $A^\top w(i)= \lambda_i w(i)$.  Write
$w=(w_1,\ldots,w_{N-p})$, we consider the function
$$F(z, w)=z^\beta [w\,w(1)]^{\gamma_1}\cdots[w\,w(N-p)]^{\gamma_{N-p}}, \qquad (z,
w)\in  \mathbb{C}^p\times\mathbb{C}^{N-p}\cap B_N,$$ where
$\beta=(\beta_1,\ldots \beta_p)\in \mathbb{N}^p$ and
$\gamma=(\gamma_1,\ldots,\gamma_{N-p})\in \mathbb{N}^{N-p}$.  It is
clear  that $[w\,w(i)]^{\gamma_i}$ is a homogeneous polynomial of
$w$ with
 degree $\gamma_i$. So that $F$ belongs to $H^2(B_N)$ and {\setlength\arraycolsep{2pt}
\begin{eqnarray*}F\circ \psi(z, w)&=&(Uz)^\beta [wA^\top w(1)]^{\gamma_1}\cdots[wA^\top
w(N-p)]^{\gamma_{N-p}}   \\ &=& e^{i(\beta_1
\theta_1+\cdots\beta_p\theta_p)}\lambda_1^{\gamma_1}\cdots\lambda_{N-p}^{\gamma_{N-p}}z^\beta
[w\,w(1)]^{\gamma_1}\cdots[w\,w(N-p)]^{\gamma_{N-p}}\\ &=&
e^{i(\beta_1
\theta_1+\cdots\beta_p\theta_p)}\lambda_1^{\gamma_1}\cdots\lambda_{N-p}^{\gamma_{N-p}}
F(z, w).\end{eqnarray*}}Thus, for any  multi-index $\alpha=(\beta,
\gamma)\in \mathbb{N}^N$, $e^{i(\beta_1
\theta_1+\cdots\beta_p\theta_p)}\lambda_1^{\gamma_1}\cdots\lambda_{N-p}^{\gamma_{N-p}}$
is an eigenvalue of $C_\psi$ with the eigenvector $F$. Hence,
$\sigma(C_\psi)$ contains the set $\bigcup\limits_{\alpha\in
\mathbb{N}^N}\Lambda_\varphi^\alpha\cup\{0\}$. If $U$ has an
irrational eigenvalue $e^{i\theta_j}$,   then the set
$\{e^{ik\theta_j}\lambda_1^{\gamma_1}\cdots\lambda_{N-p}^{\gamma_{N-p}}:
k\in \mathbb{N}\}$ is dense in the circle
$|\lambda_1|^{\gamma_1}\cdots|\lambda_{N-p}|^{\gamma_{N-p}}$. Note
that this set is contained in  $\sigma(C_\psi)$  and
$\sigma(C_\psi)$  is closed, it follows that $\sigma(C_\psi)$
contains $\bigcup\limits_{\alpha\in
\mathbb{N}^N}\Lambda_\varphi^\alpha\mathbb{T}\cup\{0\}$.

For another direction,  as in the proof of Lemma 5.4 in \cite{B10},
we need the following notations. For $\gamma\in \mathbb{N}^{N-p}$,
let $H_\gamma$ be the set of all functions $F$ in $H^2(B_N)$, which
have the form $F_\gamma(z)w^\gamma$. It is clear that
$F_\gamma(z)w^\gamma$ and $F_\beta(z)w^\beta$ are orthogonal if
$\gamma\ne \beta\in \mathbb{N}^{N-p}$. So we have the orthogonal
decomposition $H^2(B_N)=\bigoplus^{\bot}_{\gamma\in
\mathbb{N}^{N-p}}H_\gamma$. If we  set
$K_n=\bigoplus^{\bot}_{|\gamma|\ge n}H_\gamma$ for  any $n\ge 0$,
then $H^2(B_N)=\bigoplus_{|\gamma|<n}H_\gamma\oplus K_n$. Write
$F=F_\gamma(z)w^\gamma$, we see that
$$F\circ\psi(z,
w)=F_\gamma(Uz)\prod\limits_{j=1}^{N-p}\biggl(\lambda_jw_j+\sum\limits_{k>j}a_{j,k}w_k\biggr)^{\gamma_j}.$$
After introducing a natural order on the set $\mathbb{N}^{N-p}$ (see
[11, p.26] or \cite{B10}), if we use  this ordering for the
decomposition $H^2(B_N)=\bigoplus_{|\gamma|<n}H_\gamma\oplus K_n$,
then the matrix of $C_\psi$ is upper-triangular. Set $\rho(z,
w)=(Uz, w)$, let $T_\gamma$ and $S_\gamma$ respectively  denote the
diagonal blocks of $C_\psi$  and $C_{\rho}$ corresponding to
$H_\gamma$, then  we have
$T_\gamma=\lambda_1^{\gamma_1}\cdots\lambda_{N-p}^{\gamma_{N-p}}S_\gamma$.
Since $\rho$ is an elliptic automorphism of $B_N$, by Theorem 7.6 of
\cite{CM1}, the spectrum of $C_\rho$ on $H^2(B_N)$ is the closure of
all possible products of the eigenvalues of $U$, which denoted by
$X$. If $U$ has an irrational eigenvalue of modulus $1$, then
$X=\mathbb{T}$.  Applying Lemma 7.17 of \cite{CM1} or Lemma 5.3 of
\cite{B10}, we get $\sigma(S_\gamma)\subset X$ and
{\setlength\arraycolsep{2pt}
\begin{eqnarray*}\sigma(C_\psi)& \subset &\bigcup\limits_{|\gamma|<n}
\sigma(T_\gamma)\cup \sigma(C_{\psi |K_n}) \\ &\subset
&\bigcup\limits_{|\gamma|<n} \lambda_1^{\gamma_1}\cdots
\lambda_{N-p}^{\gamma_{N-p}}X\cup \sigma(C_{\psi |K_n}).
\end{eqnarray*}}Since the spectrum of $C_{\psi |K_n}$ must be
contained in the disk of radius $||C_{\psi |K_n}||$, we will obtain
the desired result if we can show that $||C_{\psi |K_n}||$ tends to
zero as $n$ tends to infinity.

Let $A=T\Sigma V$ be a singular value decomposition of $A$, where
$T, V$ are unitary and $\Sigma$ is diagonal with the diagonal
entries $\mu_1, \ldots, \mu_{N-p}$.  We see that $\mu_1, \ldots,
\mu_{N-p}$  are the non-negative square roots of the eigenvalues of
$AA^\ast$ and $\mu=\max\{\mu_1, \ldots, \mu_{N-p}\}<1$. Set
$\phi_1(z, w)=(Uz, Vw)$, $\psi_\Sigma(z, w)=(z, \Sigma w)$ and
$\phi_2(z, w)=(z, Tw)$,  so that
$C_\psi=C_{\phi_1}C_{\psi_\Sigma}C_{\phi_2}$. Note that $\phi_1$ and
$\phi_2$ are automorphisms of $B_N$, it follows that $C_{\phi_1}$
and $C_{\phi_2}$ are invertible on $H^2(B_N)$. on the other hand,
$K_n$ and $K_n^\bot$ are stabled by  $C_{\phi_1}$ and $C_{\phi_2}$.
Therefore, we only need to calculate  that $||C_{\psi_\Sigma
|K_n}||$ goes to zero as $n\to \infty$.

Let $F=\sum\limits_{|\gamma|\ge n}  F_\gamma(z)w^\gamma\in K_n$.
Composing it by $\psi_\Sigma $, we get $$F\circ \psi_\Sigma (z,
w)=\sum\limits_{|\gamma|\ge n}F_\gamma(z)(\Sigma
w)^\gamma=\sum\limits_{|\gamma|\ge n}
\mu_1^{\gamma_1}\cdots\mu_{N-p}^{\gamma_{N-p}}
F_\gamma(z)w^\gamma.$$ If we denote $\zeta=(\zeta', \zeta'')\in
\mathbb{C}^p\times\mathbb{C}^{N-p}\cap \partial B_N$, then
{\setlength\arraycolsep{2pt}
\begin{eqnarray*} ||C_{\psi_\Sigma} F||^2 &=&
\sup\limits_{0<r<1}\int_{\partial B_N} |F\circ
\psi_\Sigma(r\zeta)|^2 d\sigma(\zeta) \\ &=& \sup\limits_{0<r<1}
\sum\limits_{|\gamma|\ge n} \int_{\partial B_N}|r^{|\gamma|}
\mu_1^{\gamma_1}\cdots\mu_{N-p}^{\gamma_{N-p}}F_\gamma(r\zeta')
(\zeta'')^\gamma|^2 d\sigma(\zeta) \\ &\le &
\sum\limits_{|\gamma|\ge n}\mu^{|\gamma|}  \sup\limits_{0<r<1}
\int_{\partial B_N}|F_\gamma(r\zeta') (r\zeta'')^\gamma|^2
d\sigma(\zeta)\le \mu^n ||F||^2.
\end{eqnarray*}}Here and other places in this paper, we let  $||\cdot||$ denote the norm of $H^2(B_N)$. Applying $\mu<1$ to the above inequality
and let $n$ go to infinity, we get the desired conclusion.  \ \
$\Box$

In the second case, we find that $C_\varphi^n$ is a  compact
 operator for some positive integer $n$, so that the  spectrum of
$C_\varphi$ is easily     known.
\\ \\ {\bf Theorem 2.3.}\begin{em}  Let $\varphi\in {\mbox LFM}(B_N)$ be
elliptic, non-automorphism, and let $p=\mbox{dim}\ L_U(\varphi,
z_0)$ for its fixed point $z_0\in B_N$. If  $p=0$  and $\varphi$
fixes no boundary point. Then the spectrum of $C_\varphi$ on
$H^2(B_N)$ is the set consisting of $0$, $1$ and all possible
products of the eigenvalues of $\varphi'(z_0)$.\end{em}
\\ \par {\bf Proof.} Applying  Theorem 2.1, $\varphi$ is conjugated to a
linear fractional map $\psi$ of $B_N$  with
$$\sigma\circ \psi=A\sigma,$$ where $\sigma$ is an one-to-one
holomorphic map from $B_N$ onto   a complex ellipsoid  $\Delta_1$
and $A$ is a matrix of $\mathbb{C}^{N\times N}$ with $||A||<1$.
Moreover, in the proof of Theorem 2.1, we have $\psi(0)=0$ and
$\sigma(0)=0$.

Let $\psi^n$ denote the $n$-th iterate of $\psi$, then
$\psi^n=\sigma^{-1} A^n \sigma$ and $\psi^n(0)=0$. We claim that
there exists a positive integer $n$  so that   $C_{\psi^n}$ is
compact. Since the set $\overline{\Delta_1}=\sigma(\overline{B_N})$
is compact in $\mathbb{C}^N$ and $||A||<1$, we see that $\sup\{|A^n
\sigma(z)|: z\in \overline{B_N}\}$ goes to zero as $n$ tends to
infinity. Thus, there exist a positive integer  $M$ and a positive
constant $r<1$ such that
$$\sup\limits_{z\in \overline{B_N}}|\sigma^{-1}A^n \sigma(z)|\le r$$
holds for all $n\ge M$. That is $||\psi^n||_\infty\le r<1$ for such
$n$. As we know, a linear fractional map  $\phi$ of $B_N$ is compact
on $H^2(B_N)$ if and only if $||\phi||_\infty<1$. It follows that
$C_{\psi^n}$ is compact. Thus, by  Theorem 7.2  of \cite{CM1}, the
spectrum of  $C_{\psi^n}$ is the set consisting of $0$, $1$ and all
possible products of the eigenvalues of $(\psi^n)'(0)=\psi'(0)^n$.
On the other hand, the spectral mapping theorem gives
$[\sigma(C_\psi)]^n=\sigma(C^n_{\psi})=\sigma(C_{\psi^n})$. Hence,
the spectrum of $C_\psi$ has the desired structure.  \ \ $\Box$

For the last case in  Theorem 2.1, the structure of the spectrum of
$C_\varphi$ on $H^2(B_N)$ is similar to that in one variable with
 the symbol  fixing  a boundary  point (see
Theorem A). So we first need to compute the essential spectral
radius of $C_\varphi$. For $\zeta\in
\partial B_N$, we will use the notation $$d_\varphi(\zeta)=\liminf\limits_{z\to
\zeta}\frac{1-|\varphi(z)|}{1-|z|}$$ for an analytic map $\varphi$
of $B_N$ into itself.
\\ \\ {\bf Lemma  2.4.}\begin{em}  Let  $\varphi\in {\mbox
LFM}(B_N)$ be elliptic, non-automorphism, and let $p=\mbox{dim}\
L_U(\varphi, z_0)$ for its  fixed point $z_0\in B_N$. If  $p=0$  and
$\varphi$ fixes only one boundary point. Then the essential norm of
$C_\varphi$ on $H^2(B_N)$ satisfies  the following inequalities:
$$\limsup\limits_{|z|\to
1}\biggl(\frac{1-|z|^2}{1-|\varphi(z)|^2}\biggr)^{\frac{N}{2}}\le
||C_\varphi||_e\le C \limsup\limits_{|z|\to
1}\biggl(\frac{1-|z|^2}{1-|\varphi(z)|^2}\biggr)^{\frac{N}{2}}$$ for
a positive constant $C$. So the essential spectral radius of
$C_\varphi$ is
$$r_e(C_\varphi)=\lim\limits_{n\to\infty}\biggl(\limsup\limits_{|z|\to
1}\biggl(\frac{1-|z|^2}{1-|\varphi^n(z)|^2}\biggr)^{N/2}\biggr)^{\frac{1}{n}}.$$\end{em}
\par {\bf Proof.} Conjugation by a unitary map, we may assume that
$\varphi$ fixes $0$ and the boundary point $e_1$. We first claim
that $|\varphi(\zeta)|<1$ for any boundary point $\zeta\ne e_1$.
Otherwise, if there exist two boundary points $\zeta\ne e_1$
 and $\eta$ such that $\varphi(\zeta)=\eta$, we will get a
contradiction.

Let  $L(e_1, \zeta)=\{ce_1+(1-c)\zeta: c\in \mathbb{C}\}\cap B_N$
denote the one-dimensional affine subset of $B_N$ determined by
$e_1$ and $\zeta$  and let $[e_1]$ denote the slice  in $B_N$
through $0$ and $e_1$. If $\eta\ne e_1$, then $\varphi$ maps the
slice $L(e_1, \zeta)$ onto $L(e_1, \eta)$ with a boundary fixed
point $e_1$.  When $\zeta$ is on the slice $[e_1]$, since
$\varphi(0)=0$, we see that  $L(e_1, \zeta)=L(e_1, \eta)=[e_1]$.
Thus, $\varphi$ restricting to $[e_1]$ must be identity, which
contradicts to  $p=0$. Otherwise,  assume that $\tau_1$ and $\tau_2$
are
 automorphisms of $B_N$ fixing  $e_1$, and they map $[e_1]$ onto $L(e_1,
\zeta)$  and $L(e_1, \eta)$ onto $[e_1]$ respectively, such that the
restriction of $\tau_2\circ \varphi\circ \tau_1$ to the slice
$[e_1]$ is a linear fractional self-map of this disk onto itself
with only one boundary fixed point $e_1$. Moreover, its boundary
boundary dilation coefficient at $e_1$ satisfies $\alpha\ge
d_{\varphi}(e_1)$. Note that $\varphi$ is an elliptic linear
fractional map of $B_N$ with $\varphi(0)=0$, which  gives
$d_{\varphi}(e_1)>1$. Thus,  we have  $\alpha>1$. It is impossible
for an automorphism of the  disk with only one  boundary fixed
point.

If  $\eta=e_1$, we see that $\varphi(e_1)=\varphi(\zeta)=e_1$  and
$\varphi$ maps $L(e_1, \zeta)$ into a slice $S$   in $B_N$ through
$e_1$. Then the restriction of $\tau_2\circ \varphi\circ\tau_1$ to
$[e_1]$ is a linear fractional self-map of this disk, where $\tau_1$
and $\tau_2$ are automorphisms of $B_N$ mapping $[e_1]$ onto $L(e_1,
\zeta)$  and $S$ onto $[e_1]$ with $\tau_1(e_1)=\tau_2(e_1)=e_1$.
Suppose $\tau_1(\alpha)=\zeta$, then $\tau_2\circ
\varphi\circ\tau_1|_{[e_1]}$ is a linear factional map of the  disk
mapping  two boundary points $e_1$ and $\alpha$ to $e_1$. So it  is
a constant with  modulus $1$. Therefore, $|\tau_2\circ
\varphi\circ\tau_1(z)|=1$ for any point $z\in [e_1]$, that is,
$|\varphi(z)|=1$ for any $z\in L(e_1, \zeta)$, which contradicts to
the fact $\varphi(B_N)\subset B_N$. Thus, we show the claim.

Next, we will estimate the essential norm of $C_\varphi$. Since
$||C_\varphi||_e=||C^\ast_\varphi||_e=\inf\{||C^\ast_\varphi-F||: F\
\mbox{is compact}\}$, it is easy to see
that{\setlength\arraycolsep{2pt}
\begin{eqnarray*} ||C_\varphi||_e&\ge & \limsup\limits_{|z|\to
1}\biggl|\biggl|C_\varphi^\ast\frac{K_z}{||K_z||}\biggr|\biggr|=\limsup\limits_{|z|\to
1}\frac{||K_{\varphi(z)}||}{||K_z||}\\ &=&\limsup\limits_{|z|\to
1}\biggl(\frac{1-|z|^2}{1-|\varphi(z)|^2}\biggr)^{N/2},
\end{eqnarray*}}where $||T||_e$ denotes the essential norm of a operator $T$ and $K_z(w)=\frac{1}{(1-<z, w>)^N}$ is the
reproducing kernel of $H^2(B_N)$.

We introduce the pullback measure $\tau_\varphi$ on $\overline{B_N}$
defined by $\tau_\varphi(E)=\sigma[(\varphi^\ast)^{-1}(E)]$ for
Borel subsets $E$ of $\overline{B_N}$,  where $\varphi^\ast$ denotes
the radial limit of $\varphi$. Then $$\int_{\partial B_N} f\circ
\varphi^\ast (\zeta)d \sigma (\zeta)=\int_{\overline{B_N}} f d
\tau_\varphi$$ for every Borel function $f\ge 0$ on
$\overline{B_N}$. Let $I_{\tau_\varphi}$ denote the densely defined
inclusion operator of $H^2(B_N)$ into $L^2(\tau_\varphi)$. Using
Theorem 5.1  of \cite{Choe}, we get
$$||C_\varphi||_e^2=||I_{\tau_\varphi}||^2_e\cong
||\tau_\varphi||_{e,C},$$ where the ``essential Carleson norm" is
defined by $$||\tau_\varphi||_{e,C}=\limsup\limits_{t\to
0}\sup\limits_{\zeta\in \partial
B_N}\frac{\tau_\varphi(\Omega(\zeta, t))}{\sigma(Q(\zeta, t))}$$
with $\Omega(\zeta, t)=\{z\in\overline{B_N}: |1-<z, \zeta>|<t\}$ and
$Q(\zeta, t)=\Omega(\zeta, t)\cap \partial B_N$.

Since $\varphi(e_1)=e_1$ and $\varphi$ maps any boundary point
$\zeta\ne e_1$  into $B_N$, this gives
$$||\tau_\varphi||_{e,C}=\limsup\limits_{t\to
0}\frac{\tau_\varphi(\Omega(e_1, t))}{\sigma(Q(e_1,
t))}=\limsup\limits_{t\to
0}\frac{\sigma[(\varphi^\ast)^{-1}(\Omega(e_1, t)]}{\sigma(Q(e_1,
t))}.$$ Applying  similar arguments as for proving $\varphi(Q(\zeta,
t))\subset\Omega(\eta, At)$ with $\varphi(\zeta)=\eta$ in Lemma 3.40
of \cite{CM1}, where $A$  is a constant depending on
$d_\varphi(\zeta)$, we can calculate that
$$||\tau_\varphi||_{e,C}=\limsup\limits_{t\to
0}\frac{\sigma[(\varphi^\ast)^{-1}(\Omega(e_1, t)]}{\sigma(Q(e_1,
t))}=\limsup\limits_{z\to
e_1}\biggl(\frac{1-|z|^2}{1-|\varphi(z)|^2}\biggr)^N$$ because of
$\sigma(Q(e_1, t))\sim t^N$.  Therefore, $$||C_\varphi||^2_e\le
C||\tau_\varphi||_{e, C}\le C\limsup\limits_{|z|\to
1}\biggl(\frac{1-|z|^2}{1-|\varphi(z)|^2}\biggr)^N$$
 holds for a positive constant $C$. Note that the essential spectral
 radius of $C_\varphi$ can be computed by
$r_e(C_\varphi)=\lim\limits_{n\to\infty}{||C^n_\varphi||_e^{1/n}}=\lim\limits_{n\to\infty}{||C_{\varphi^n}||_e^{1/n}}$
and we have obtained  that $$||C_{\varphi^n}||_e\cong
\limsup\limits_{|z|\to
1}\biggl(\frac{1-|z|^2}{1-|\varphi^n(z)|^2}\biggr)^{N/2}.$$ Thus, we
get
$$r_e(C_\varphi)=\lim\limits_{n\to\infty}{||C_{\varphi^n}||_e^{1/n}}=\lim\limits_{n\to\infty}\biggl(\limsup\limits_{|z|\to
1}\biggl(\frac{1-|z|^2}{1-|\varphi^n(z)|^2}\biggr)^{N/2}\biggr)^{\frac{1}{n}}.$$
$\Box$

Kamowitz  \cite{Kam} was the first to investigate the spectrum of
$C_\varphi$ on $H^2(D)$ when $\varphi$ has an interior  fixed point
$a$, he proved that $\sigma(C_\varphi)=\{\lambda: |\lambda|\le
r_e(C_\varphi)\}\cup\{\varphi'(a)^n: n=1,2,\ldots\}\cup\{1\}$ if
$\varphi$ is analytic in a neighborhood of $D$, not an inner
function. In \cite{CM94} Cowen and MacCluer obtained the same
spectral structure for $C_\varphi$ when $\varphi$ is univalent, not
automorphism, and  $\varphi(a)=a\in B_N$. From their results one can
easily deduce the spectrum of $C_\varphi$ if $\varphi$ is   a
 linear fractional map of $D$   with  an interior and a  boundary fixed
points (see Theorem A). Moreover, the method of Cowen and MacCluer
was used by many authors to study the spectrum of $C_\varphi$ on
other spaces for $\varphi$ in the same case, For example, on
$H^\infty(D)$, the Bloch space and BMOA in one variable, see
\cite{Zheng}, \cite{MS},  \cite{LP}; for a generalization on
$H^\infty(B_E)$ and $H^\infty(B_N)$ see \cite{GGL} and \cite{YZ},
where $B_E$ is an open unit ball on  a complex Banach space. We will
use same ideas and approaches suggested by the work of Kamowitz and
that of Cowen and MacCluer to obtain the following result.
\\ \\ {\bf Theorem 2.5.}\begin{em}  Let $\varphi$ be the same as in Lemma 2.4.
Then the spectrum of $C_\varphi$ on $H^2(B_N)$ is
$$\sigma(C_\varphi)=\{\lambda: |\lambda|\le \rho\}\cup\{\mbox{all possible products of the eigenvalues of}\ \varphi'(z_0)\}\cup\{1\},$$
where $\rho$ is the essential spectral radius of
$C_\varphi$.\end{em}
\\ \par
As in the proof of Lemma 2.4, we may assume $\varphi(0)=0$ and
$\varphi(e_1)=e_1$. For non-negative integers $m$, let $\mathcal
{H}_m$ be the subspace of $H^2(B_N)$ spanned by the monomials of
total degree greater than or equal to $m$,  that is,  any function
$f\in \mathcal {H}_m$ can be written as $f=\sum_{|s|\ge m} f_s$,
where $f_s$ is a homogeneous polynomial of degree $s$. Obviously,
$\mathcal {H}_m$ is invariant for $C_\varphi$. A easy computation
shows that the reproducing kernel in $\mathcal {H}_m$ at the point
$w$ is
$$K^m_w(z)=\sum\limits_{|s|\ge m}\frac{(N-1+s)!}{(N-1)!s!}<z, w>^s.$$
In particular, we have $$||K^m_w||^2=\sum\limits_{|s|\ge
m}\frac{(N-1+s)!}{(N-1)!s!}|w|^{2s}$$ and $$C'_{N,m} |w|^m||K_w||\le
||K^m_w||\le C''_{N,m} |w|^m||K_w||,$$  where $C'_{N,m}$ and
$C''_{N,m}$ are positive constants only depending on $N$ and $m$.

Just as in the disk, our generalization for Proposition 7.32 of
\cite{CM1} to the unit ball is the following lemma. We only need a
similar argument used in computing the spectrum of a compact
composition operator on $H^2(B_N)$, so the proof will be omitted.
\\ \\ {\bf Lemma 2.6.}\begin{em} For $\varphi$ as in Lemma 2.4, the spectrum of
$C_\varphi$ on $H^2(B_N)$ contains $1$ and all possible products of
the eigenvalues of $\varphi'(z_0)$. Moreover, if $\lambda\ne 0$ is
an eigenvalue of $C_\varphi$, then $\lambda$ will be a product of
the eigenvalues of  $\varphi'(z_0)$.\end{em}
\\ \par
We say the sequence of points $\{z_k\}^\infty_{-K}$ is an iteration
sequence for $\varphi$ if $\varphi(z_k)=z_{k+1}$ for $k\ge -K$. In
the proof of Theorem 2.5, the following fact will be needed.
\\ \\ {\bf Lemma D. \cite{CM94}}\begin{em}  If $\varphi$ maps the unit ball into
itself, $\varphi(0)=0$  and $\varphi$ is not unitary on any slice in
$B_N$. Suppose $r<1$ is given and $\{z_k\}^\infty_{-K}$ is any
iteration sequence with $|z_n|\le r$. Then there exists $c<1$ such
that $\frac{|z_{k+1}|}{|z_k|}\le c$ for all $k\ge n$.\end{em}
\\ \\ {\bf Lemma E. \cite{CM94}}\begin{em}
 For $\varphi$  in Lemma D  and
$0<r<1$. There exists $M<\infty$ such that if $\{z_k\}^\infty_{-K}$
is an iteration sequence with $|z_l|\ge r$ for some $l\ge 0$ and if
$\{w_k\}^\infty_{-K}$ is arbitrary, then there is $h$ in
$H^\infty(B_N)$  such that $h(z_k)=w_k$ for  $-K\le k\le l$ and
$||h||_\infty\le M\sup\{|w_k|:-K\le k\le l\}$. \end{em}
\\ \par {\bf
Proof of Theorem 2.5.} Using  Lemma 2.6 and the fact that the
spectrum is closed, it suffices to prove $$\{\lambda: |\lambda|\le
\rho\}\subset \sigma(C_\varphi).$$ Let $C_m=C_{\varphi | \mathcal
{H}_m}$ and assume $0<|\lambda|<\rho$. By Lemma 7.17 of \cite{CM1},
if we can show that $\lambda$ is in the spectrum of $C_m$ for some
positive integer $m$, then we get the desired conclusion. So we will
try to find a positive integer $m$ such that $C_m^\ast-\lambda I$ is
not bounded from below.

Fixing $\delta$ with $0<\delta<1$, suppose we have an iteration
sequence $\{z_k\}^\infty_{-K}$ with $|z_0|>\delta$. Let $n=\max\{k:
|z_k|\ge \delta\}$. By Lemma D, we can choose $c$ with
$\sqrt{\delta}<c<1$ such that $|z_{k+1}|<c|z_k|$ for
$|z_k|<\sqrt{\delta}$. On the other hand, if
$|z_n|>\sqrt{\delta}>\delta$, then
$|z_{n+1}|<\delta<\sqrt{\delta}|z_n|<c|z_n|$. Thus, we have
$|z_{k+1}|<c|z_k|$ for all $k\ge n$.  By induction,
$|z_k|<c^{k-n}|z_n|$ holds for all $k\ge n$.

Now, let's define
$L_\lambda=\sum\limits^\infty_{k=0}\lambda^{-(k+1)}K_{z_k}^m$, we
will see that $L_\lambda$ is well defined and it is bounded. For
$k>n$, $|z_k|<\delta$ gives that $$||K_{z_k}^m||\le  C''_{N,m}
|z_k|^m||K_{z_k}||=C''_{N,m} \frac{|z_k|^m}{(1-|z_k|)^{N/2}}\le C
|z_k|^m,$$ where $C$ is a positive constant depending on $N, m$ and
$\delta$. It follows that
$$\sum\limits^\infty_{k=n+1}|\lambda|^{-(k+1)}||K_{z_k}^m||\le C\sum\limits^\infty_{k=n+1}\frac{|z_k|^m}{|\lambda|^{k+1}}
\le
C\frac{|z_n|^m}{|\lambda|^{n+1}}\sum\limits^\infty_{k=n+1}\biggl(\frac{c^m}{|\lambda|}\biggr)^{k-n}.$$
Since $c<1$, if we choose $m$ so large that $c^m<|\lambda|$, then
the series defining $L_\lambda$  converges.

Next, we will estimate $||(C_m^\ast-\lambda I)
L_\lambda||/||L_\lambda||$. Since $\mathcal {H}_m$ is invariant for
$C_\varphi$, it is easy to see  $C_\varphi^\ast
K_w^m=K_{\varphi(w)}^m$ and  so {\setlength\arraycolsep{2pt}
\begin{eqnarray*}(C_m^\ast-\lambda I)
L_\lambda &=&(C_m^\ast-\lambda I)
\sum\limits^\infty_{k=0}\lambda^{-(k+1)}K_{z_k}^m=\sum\limits^\infty_{k=0}\biggl(\lambda^{-(k+1)}K_{\varphi(z_k)}^m-\lambda^{-k}K_{z_k}^m\biggr)
\\ &=&
\sum\limits^\infty_{k=0}\biggl(\lambda^{-(k+1)}K_{z_{k+1}}^m-\lambda^{-k}K_{z_k}^m\biggr)=-K_{z_0}^m.
\end{eqnarray*}}This means  $||(C_m^\ast-\lambda I)
L_\lambda||=||K_{z_0}^m||$. It remains to give a lower bound for
$||L_\lambda||$.

Choose an $m$-homogenous polynomial $P(z)$ on $B_N$ satisfying
$||P||_\infty=1$ and $|P(z_n)|=|z_n|^m$.  It is clear that
$|P(z)/|z|^m|=|P(z/|z|)|\le 1$  for any $z\in B_N$, which gives
$|P(z)|\le |z|^m$. Applying Lemma E to the iteration sequence
$\{z_k\}^\infty_{-K}$, there is a function $g\in H^\infty(B_N)$ such
that $||g||_\infty\le M$,  $g(z_k)=0$  for $0\le k<n$,  and
$g(z_n)=1$. For $f=\frac{Pg}{(1-<z, z_n>)^N}$, we obtain
$$<L_\lambda,
f>=\biggl<\sum\limits^\infty_{k=0}\lambda^{-(k+1)}K_{z_k}^m,f\biggr>
=\lambda^{-(n+1)}\overline{f(z_n)}+\sum\limits^\infty_{k=n+1}\lambda^{-(k+1)}\overline{f(z_k)}.$$
Note that $|z_k|< \delta$ for $k>n$, so
$|1-<z_k,z_n>|^N\ge(1-|z_k|)^N\ge C(\delta)$. Choose $m$ large
enough so that $c^m/|\lambda|<\frac{1}{2^{N+2}M}$, we calculate that
{\setlength\arraycolsep{2pt}
\begin{eqnarray*}\biggl|\sum\limits^\infty_{k=n+1}\lambda^{-(k+1)}\overline{f(z_k)}\biggr|&=&
\biggl|\sum\limits^\infty_{k=n+1}\frac{1}{\lambda^{k+1}}\overline{\frac{P(z_k)g(z_k)}{(1-<z_k,z_n>)^N}}\biggr|
\le\sum\limits^\infty_{k=n+1}\frac{|P(z_k)g(z_k)|}{|\lambda|^{k+1}|1-<z_k,z_n>|^N}\\
&\le &
\sum\limits^\infty_{k=n+1}\frac{M|z_k|^m}{|\lambda|^{k+1}|1-<z_k,z_n>|^N}\le
\frac{M|z_n|^m}{C(\delta)|\lambda|^{n+1}}
\sum\limits^\infty_{k=n+1}\biggl(\frac{c^m}{|\lambda|}\biggr)^{k-n}
\\ &\le& \frac{|z_n|^m}{2^{N+1}C(\delta)|\lambda|^{n+1}}.
\end{eqnarray*}}Therefore,{\setlength\arraycolsep{2pt}
\begin{eqnarray*}|<L_\lambda,
f>|&\ge &
\biggl|\lambda^{-(n+1)}f(z_n)\biggr|-\biggl|\sum\limits^\infty_{k=n+1}\lambda^{-(k+1)}\overline{f(z_k)}\biggr|
\\ &\ge &
\frac{1}{2}\frac{|P(z_n)g(z_n)|}{|\lambda|^{n+1}(1-|z_n|^2)^N}+
\frac{1}{2}\frac{|P(z_n)g(z_n)|}{|\lambda|^{n+1}(1-|z_n|^2)^N}-\frac{|z_n|^m}{2^{N+1}C(\delta)|\lambda|^{n+1}}
\\ &= &
\frac{|z_n|^m}{2|\lambda|^{n+1}}\frac{1}{(1-|z_n|^2)^N}+\frac{|z_n|^m}{2|\lambda|^{n+1}}\frac{1}{(1-|z_n|^2)^N}-\frac{|z_n|^m}{2^{N+1}C(\delta)|\lambda|^{n+1}}
\\&\ge &\frac{|z_n|^m}{2|\lambda|^{n+1}}\frac{1}{(1-|z_n|^2)^N},
\end{eqnarray*}}where we have used the fact that $|z_n|\ge \delta$ gives
$(1-|z_n|^2)^N\le 2^NC(\delta)$ to obtain the third inequality.
Since
$$||f||=\biggl|\biggl|\frac{Pg}{(1-<z, z_n>)^N}\biggr|\biggr|\le
\frac{M}{(1-|z_n|^2)^{N/2}},$$ we deduce that
{{\setlength\arraycolsep{2pt}
\begin{eqnarray*}||L_\lambda||&\ge &\frac{|<L_\lambda,
f>|}{||f||}\ge\frac{\frac{|z_n|^m}{2|\lambda|^{n+1}}\frac{1}{(1-|z_n|^2)^N}}{\frac{M}{(1-|z_n|^2)^{N/2}}}
\\ &=&\frac{|z_n|^m}{2M|\lambda|^{n+1}}\frac{1}{(1-|z_n|^2)^{N/2}}=\frac{|z_n|^m||K_{z_n}||}{2M|\lambda|^{n+1}}
\\ &\ge &\frac{||K_{z_n}^m||}{2MC'_{N,m}|\lambda|^{n+1}}.
\end{eqnarray*}}Hence, $$\frac{||(C_m^\ast-\lambda I)
L_\lambda||}{||L_\lambda||}\le
\frac{2MC'_{N,m}|\lambda|^{n+1}||K_{z_0}^m||}{||K_{z_n}^m||}.$$ When
$m$ is fixed, we see that $||K_w^m||\le ||K_w||\le ||K_w^m||+C_m$
for all $w\in B_N$, which gives
$$\frac{||K_{\varphi^n(w)}^m||}{||K_w^m||}\approx\frac{||K_{\varphi^n(w)}||}{||K_w||}.$$
On the other hand, by Lemma 2.4, the essential spectral radius of
$C_\varphi$ is
$$\rho=\lim\limits_{n\to\infty}\biggl(\limsup\limits_{|w|\to
1}\biggl(\frac{1-|w|^2}{1-|\varphi^n(w)|^2}\biggr)^{N/2}\biggr)^{\frac{1}{n}}=\lim\limits_{n\to\infty}\biggl(\limsup\limits_{|w|\to
1}\frac{||K_{\varphi^n(w)}||}{||K_w||}\biggr)^{\frac{1}{n}}.$$ The
above arguments show that
$$\rho=\lim\limits_{n\to\infty}\biggl(\limsup\limits_{|w|\to
1}\frac{||K^m_{\varphi^n(w)}||}{||K^m_w||}\biggr)^{\frac{1}{n}}.$$
Since $\rho>0$, we   have $\limsup\limits_{|w|\to
1}||K^m_{\varphi^n(w)}||=\infty$ for every $n$. Thus, for any
positive integer $n$, there exist points $w$ near $\partial B_N$
such that $|\varphi^n(w)|>\delta$. Moreover,  there exists  $\rho'$
with $|\lambda|<\rho'<\rho$ such that
$$\frac{||K^m_{\varphi^n(w)}||}{||K^m_w||}\ge (\rho')^n.$$
Let $z_0=w$ and $z_{k+1}=\varphi(z_k)$ for $k\ge 0$. It is clear
that  the iteration sequence $\{z_k\}^\infty_0$ satisfies
$|z_0|>|z_n|>\delta$. Consequently, for this iteration sequence, we
obtain {\setlength\arraycolsep{2pt}
\begin{eqnarray*}\frac{||(C_m^\ast-\lambda I)
L_\lambda||}{||L_\lambda||}&=&\frac{||(C_m^\ast-\lambda I)
\sum^\infty_{k=0}\lambda^{-(k+1)}K_{z_k}^m||}{||\sum^\infty_{k=0}\lambda^{-(k+1)}K_{z_k}^m||}
\\ & \le & \frac{2MC'_{N,m}|\lambda|^{n+1}||K_{z_0}^m||}{||K_{z_n}^m||}\le
2MC'_{N,m}|\lambda|\biggl(\frac{|\lambda|}{\rho'}\biggr)^n.
\end{eqnarray*}}Note that $|\lambda|<\rho'$, we may form iteration sequences for which $n$ is
sufficiently large, so that the above inequality is sufficiently
small. This forces that $C_m^\ast-\lambda I$ is not bounded below.
Hence, $C_m^\ast-\lambda I$ is not invertible and so we complete the
proof. \ \ $\Box$

\section*{\large 3 Spectra of hyperbolic  composition operators}

Let  $\varphi$ be a hyperbolic linear fractional self-map  of $B_N$,
it may have one or two boundary fixed points. This implies that the
associated composition operator $C_\varphi$ may have different
characterization for its spectrum on $H^2(B_N)$. First, we will find
a conjugated form of $\varphi$ on the Siegel half-plane for each
case.

Recall that the unit ball $B_N$ is biholomorphic to the Siegel
half-plane $H_N=\{(z, w)\in \mathbb{C}\times\mathbb{C}^{N-1}:
\mbox{Re}\, z>|w|^2\}$ via the Cayley transform $\sigma_C$  defined
by $$\sigma_C(z, w)=\biggl(\frac{1+z}{1-z}, \frac{w}{1-z}\biggr),
\qquad (z, w)\in \mathbb{C}\times\mathbb{C}^{N-1}.$$ It extends to a
homeomorphism of $\overline{B_N}$ onto $H_N\cup
\partial H_N\cup \{\infty\}$, the one-point compactification of
$\overline{H_N}$. Its reciprocal is given by $$\sigma_C^{-1}(z,
w)=\biggl(\frac{z-1}{z+1}, \frac{2w}{z+1}\biggr).$$

Let $\varphi$ be a linear fractional map of $B_N$ with Denjoy-Wollf
point $e_1$ and  boundary dilatation  coefficient $\alpha$. Applying
Proposition 4.2 of \cite{BCD}, $\varphi$ is conjugated to a map
$\psi$ on $H_N$ with the form $$\psi(z, w)=\frac{1}{\alpha}(z+<w,
b>+c, Aw+d), \qquad (z, w)\in H_N,$$ where $c\in\mathbb{C}$,
$b,d\in\mathbb{C}^{N-1}$ and $A\in \mathbb{C}^{(N-1)\times(N-1)}$,
which  satisfy $\alpha \mbox{Re}\, c\ge |d|^2$ and
$||A||\le\sqrt{\alpha}$. If $\varphi$ is hyperbolic, then
$\alpha<1$. Combining this with
$||A^\ast||=||A||\le\sqrt{\alpha}<1$, we see that $A^\ast-I$ is
invertible. Thus, there exists a vector $k_1\in \mathbb{C}^{N-1}$
such that $b=2(A^\ast-I)k_1$. Consider the following Heisenberg
transformation $$\eta(z, w)=(z+2<w, k_1>+k_2, w+k_1),\qquad (z,
w)\in H_N,$$ with $\mbox{Re}\, k_2=|k_1|^2$, it is an automorphism
of $H_N$. Conjugating  $\psi$ by $\eta$, we obtain
$$\phi(z, w)=\eta^{-1}\circ\psi\circ\eta(z, w)=\frac{1}{\alpha}(z+c', Aw+d'),
\qquad (z, w)\in H_N,$$ where  $c'\in \mathbb{C}$,
$d'\in\mathbb{C}^{N-1}$ and $\alpha\mbox{Re}\, c'\ge |d'|^2.$ If
$c'=c_1+i c_2$ with $c_1, c_2\in \mathbb{R}$, let $\nu(z,
w)=(z-\frac{ic_2}{1-\alpha}, w)$,  then $\nu$ is an automorphism of
$H_N$ and $$\nu^{-1}\circ\phi\circ\nu(z, w)=\frac{1}{\alpha}(z+c_1,
Aw+d'),\qquad (z, w)\in H_N,$$ with $\alpha c_1 \ge |d'|^2$.   Thus,
we may assume $c'\in \mathbb{R}$.

Therefore, similar to the classification of semi-group of hyperbolic
linear fractional self-maps of $B_N$ in \cite{BCD}, we have the
following result.
\\ \\ {\bf Proposition 3.1.}\begin{em} Let $\varphi$ is a hyperbolic linear
fractional self-map of $B_N$ with the boundary dilatation
coefficient $\alpha$.

(i) If $\varphi$ fixes only one boundary point, then $\varphi$ is
conjugated to a map $\phi$ of the form$$\phi(z,
w)=\frac{1}{\alpha}(z+c, Aw+d),\qquad (z, w)\in H_N,$$ where $c\in
\mathbb{R}$, $d\in \mathbb{C}^{N-1}$ with $\alpha c\ge |d|^2$ and
$A\in \mathbb{C}^{(N-1)\times (N-1)}$ satisfying $||A||\le
\sqrt{\alpha}$.

(ii) If $\varphi$ fixes two boundary points, then $\varphi$ is
conjugated to a map $\phi$ of the form $$\phi(z,
w)=(\frac{z}{\alpha}, \frac{Aw}{\sqrt{\alpha}}),\qquad (z, w)\in
H_N,$$ where $A$ is a matrix of order $N-1$ with $||A||\le
1$.\end{em}
\\ \par
{\bf Proof.}  The previous arguments give that $\varphi$ is
conjugated to a map on  the Siegel half-plane $H_N$ with the form
$$\phi(z, w)=\frac{1}{\alpha}(z+c, Aw+d),$$ where $c\in \mathbb{R}$ with $\alpha c\ge |d|^2$ and
$||A||\le \sqrt{\alpha}$. If $\varphi$ fixes only one boundary
point, then $c\ne 0$. To see this, we may assume $c=0$. It is clear
that $d=0$ from $|d|^2\le \alpha c$. A calculation shows that the
conjugate map $\phi$  of $\varphi$ fixes two boundary points  $0$
and $\infty$. So $\varphi$ must have two boundary fixed points  and
we get a contradiction.

If $\varphi$ has two boundary fixed points, so is its conjugate map
$\phi$. Let $(z, w)\in\partial H_N$  be another fixed point of
$\phi$ different from its Denjoy-Wollf point. We   see that
$\mbox{Re}\, z=|w|^2$, $z+c=\alpha z$ and $Aw+d=\alpha w$. Since
$c\in \mathbb{R}$  and $\alpha<1$, then $c=(\alpha-1)z$ gives $c\le
0$. On the other hand, we have $\alpha c\ge |d|^2\ge 0$. Thus,
$c=d=0$ and so the conjugate map $\phi$ of $\varphi$ has the form
$$\phi(z, w)=\frac{1}{\alpha}(z, Aw)=(\frac{z}{\alpha}, \frac{A'w}{\sqrt{\alpha}})$$ with  $||A'||=||\frac{1}{\sqrt{\alpha}}A||\le
1$.  \ \ $\Box$

Now, we will compute the spectra of hyperbolic composition operators
on $H^2(B_N)$. For the case $\varphi$ fixing only one boundary
point, we need  the technique used in determining the spectra of
composition operators on some weighted Hardy spaces, whose symbols
are hyperbolic automorphisms of $B_N$ (see Theorem 7.7 of
\cite{CM1}). Moreover, the spectrum of $C_\varphi$ has similar
properties as  that in the disk. This result has been obtained on
some weighted Hardy spaces in the first author's doctoral thesis
\cite{J08}.
\\ \\ {\bf Theorem 3.2.}\begin{em}  Let $\varphi\in\mbox{LFM} (B_N)$ be
hyperbolic, non-automorphism. If $\varphi$ has only one boundary
fixed point with the boundary dilatation coefficient $\alpha$. Then
the spectrum and the essential spectrum  of $C_\varphi$ on
$H^2(B_N)$ are the closed disk $|\lambda|\le \alpha^{-N/2}$.
Moreover, every point in the disk  $|\lambda|<\alpha^{-N/2}$ is an
eigenvalue of $C_\varphi$ of infinite multiplicity.\end{em}
\\ \par
{\bf Proof.} Proposition 3.1 tells us that $\varphi$ is conjugated
to a map  $$\phi(z, w)=\frac{1}{\alpha}(z+c, Aw+d),\qquad (z, w)\in
H_N,$$ with  $\alpha c\ge |d|^2$  and $||A||\le \sqrt{\alpha}$.
Transferring  back to the unit ball by the Cayley transform
$\sigma_C$, we have {\setlength\arraycolsep{2pt}
\begin{eqnarray*}\psi(z, w)&=& \sigma_C^{-1}\circ\phi\circ\sigma_C(z, w)
\\ &=& \biggl(\frac{1+c-\lambda+(1-c+\lambda)z}{1+c+\lambda+(1-c-\lambda)z},
\frac{2Aw+2(1-z)d}{1+c+\lambda+(1-c-\lambda)z}\biggr),\qquad (z,
w)\in B_N.
\end{eqnarray*}}

Next, we will determine the spectrum of $C_\psi$. We first prove
that each point in the disk $|\lambda|<\alpha^{-N/2}$ is an
eigenvalue of $C_\psi$ of infinite multiplicity. For $z\in D$,
define
$$\widetilde{\psi_1}(z)=\frac{1+c-\lambda+(1-c+\lambda)z}{1+c+\lambda+(1-c-\lambda)z}$$
and set
$$\widetilde{\beta}(k)^2=\frac{(N-1)!k!}{(N-1+k)!}, \qquad k=1,2,\ldots.$$
If we set $\widehat{\beta}(k)^2=(N-1)!(k+1)^{1-N}$,  then
$\widehat{\beta}(k)^2\ge \widetilde{\beta}(k)^2$ and
$H^2(\widehat{\beta}, D)\subset H^2(\widetilde{\beta}, D)$, where
$H^2(\widetilde{\beta}, D)$  and $H^2(\widehat{\beta}, D)$
respectively denote the weighted Hardy spaces with the weights
$\widetilde{\beta}$ and $\widehat{\beta}$, see Section 4 for their
definitions.

Since $c>0$, the linear fractional self-map $\widetilde{\psi_1}$ of
$D$ fixes one boundary point $1$ and the other point
$z_0=\frac{c+1-\lambda}{c-1+\lambda}$ with $|z_0|>1$.  In Theorem 8
of \cite{Hu}, Hurst has investigated the spectrum of
$C_{\widetilde{\psi_1}}$ on the weighted Hardy space $H^2(\beta, D)$
with the weight $\beta(k)=(k+1)^\alpha$ $(\alpha\le 0)$, and proved
that each point of the disk
$|\lambda|<(\widetilde{\psi_1})'(1)^{(2\alpha-1)/2}$ is an
eigenvalue of $C_{\widetilde{\psi_1}}$ of infinite multiplicity. It
is clear that  $1-N\le 0$. Therefore, using Hurst's result, if
$\lambda$ is in the disk
$|\lambda|<(\widetilde{\psi_1})'(1)^{-N/2}$, we may find infinitely
many linearly independent functions $f$ in $H^2(\widehat{\beta},
D)\subset H^2(\widetilde{\beta}, D)$ such that
$f\circ\widetilde{\psi_1}=\lambda f$.  For $(z, w)\in B_N$, define
the extension operator by $Ef(z, w)=f(z)$,  by Proposition 2.21 of
\cite{CM1}, $E$ is an isometry of $H^2(\widetilde{\beta}, D)$ into
$H^2(B_N)$. Note that $\widetilde{\psi_1}(z)=\psi_1(z, w)$ and
$(\widetilde{\psi_1})'(1)=\alpha$, where $\psi_1$ is the first
coordinate  of $\psi$. Thus, we have $Ef\in H^2(B_N)$ and
{\setlength\arraycolsep{2pt}\begin{eqnarray*}Ef\circ\psi(z,
 w)&=& Ef(\psi_1(z, w),\ldots, \psi_N(z, w))=f\circ\psi_1(z, w)
 \\ &=& f\circ\widetilde{\psi}_1(z)=\lambda f(z)=\lambda Ef(z, w).
\end{eqnarray*}}This shows that $\lambda$ is an eigenvalue of
$C_\psi$ on $H^2(B_N)$  of  infinite multiplicity.

By Theorem C, we see that the spectral radius of $C_\psi$  is
$\alpha^{-N/2}$. Hence, the spectrum of $C_\psi$ is contained in the
disk $|\lambda|\le \alpha^{-N/2}$. Since the essential spectrum is
contained in the spectrum and contains all  eigenvalues of infinite
multiplicity, we obtain  the desired conclusion. \ \ $\Box$

Finally, when a hyperbolic linear fractional map $\varphi$ fixes two
boundary points, the first author \cite{J08} has proved that each
point of the annulus $\alpha^{N/2}<|\lambda|<\alpha^{-N/2}$ is an
eigenvalue of $C_\varphi$ on $H^2(B_N)$ of infinite multiplicity,
and its spectrum  is contained in the disk $|\lambda|\le
\alpha^{-N/2}$. Applying a different approach, Jury \cite{Jury}
obtained the same  result. Recently, the spectrum of $C_\varphi$ has
been determined by Xu and Deng \cite{XD},   the main idea also comes
from Bayart \cite{B10}. By Theorem 3.1, we see that $\varphi$ is
conjugated to a map on $H_N$ with the form $\phi(z,
w)=(\frac{z}{\alpha}, \frac{Aw}{\sqrt{\alpha}})$, then the spectrum
of $C_\varphi$ is  as follows.
\\ \\ {\bf Theorem F.}\begin{em}  Let $\varphi\in\mbox{LFM} (B_N)$ be hyperbolic,
non-automorphism. If $\varphi$ fixes two boundary points with the
boundary dilation coefficient $\alpha$, then the spectrum of
$C_\varphi$ on $H^2(B_N)$ is
$$\sigma(C_\varphi)=\bigcup\limits_{\beta\in \mathbb{N}^{N-1}}\Lambda^\beta S\cup \{0\},$$
where $\Lambda^\beta$ denotes $\lambda_1^{\beta_1}\cdots
\lambda_{N-1}^{\beta_{N-1}}$ with $\lambda_1,\ldots,\lambda_{N-1}$
the eigenvalues  of $A$ and  $S$ is the annulus $\{\lambda:
\alpha^{N/2}<|\lambda|<\alpha^{-N/2}\}$.\end{em}

\section*{\large 4 Concluding remarks}

Now, the   spectra of all linear fractional composition operators on
$H^2(B_N)$ have been completely determined. Moreover,  some works
are adapted to other spaces, for example, the spectra of parabolic
composition operators have been generalized to weighted Bergman
spaces by Bayart \cite{B10}.

Let's  recall that the weighted Bergman space $A^2_\alpha(B_N)$ is
the space of all analytic functions $f$ in $B_N$ such that
$$||f||^2_{A^2_\alpha}=\int_{B_N}|f(z)|^2(1-|z|^2)^\alpha dv(z)
<\infty$$ for $\alpha>-1$. For the notational conveniences, we set
$A^2_{-1}(B_N)=H^2(B_N)$. Let $f(z)=\sum_{k=0}^\infty f_k(z)$ be the
homogeneous expansion of $f$. For $s \ge 0$, the fractional
derivative of $f$ order $s$ is defined by
$$R^s f(z)=\sum\limits_{k=1}^\infty k^sf_k(z).$$ Thus, for $\alpha\ge -1$ and $s\ge 0$, the
holomorphic Sobolev space $A^2_{\alpha, s}$ is defined as
$$A^2_{\alpha, s}(B_N)=\{f\in H(B_N): R^s f\in A^2_\alpha(B_N)\}$$
 with the norm $||f||_{A^2_{\alpha, s}}=||R^sf||_{A^2_\alpha}+|f(0)|$.
Koo and Park \cite{KP} proved that $C_\varphi$ is bounded on
$A^2_{\alpha, s}(B_N)$ if and only if $\varphi$ satisfies Wogen's
condition for $\varphi\in C^{s+4}(\overline{B_N})$. As we know, any
linear fractional self-map of $B_N$ satisfies Wogen's condition (see
\cite{CM2}). This implies that all linear fractional composition
operators are bounded on $A^2_{\alpha, s}(B_N)$.

Let $H^2(\beta, B_N)$ be  the weighted Hardy space with the weight
$\beta(k)=(k+1)^v$, where $\nu$ is a real number. Its norm is given
by
$$||f||^2_{H^2(\beta)}=\sum\limits_0^{\infty}||f_k||^2\beta(k)^2.$$  It will be
convenient to work with an equivalent norm on  $H^2(\beta, B_N)$,
the following lemma is a generalization of Lemma 1.2 in \cite{GM} to
the unit ball.
\\ \\ {\bf Lemma 4.1.}\begin{em}  Suppose that $\nu$ is a real number and $s\ge
0$ is an integer such that $s\ge\nu$. Then the norm of the weighted
Hardy space $H^2(\beta, B_N)$ with the weight $\beta(k)=(k+1)^v$ is
equivalent to the norm of  the  Sobolev space $A^2_{2s-2\nu-1,
s}$.\end{em}
\\ \par
{\bf Proof.} Let $f=\sum_{k=0}^\infty
f_k(z)=\sum_{k=0}^\infty\sum_{|\alpha|=k} a_\alpha z^\alpha$ be  in
$H^2(\beta, B_N)$. We compute that
{\setlength\arraycolsep{2pt}\begin{eqnarray*}||R^sf||^2_{A^2_c}
 &=&\int_{B_N}\biggl|\sum_{k=1}^\infty\sum\limits_{|\alpha|=k}
k^sa_\alpha z^\alpha\biggr|^2(1-|z|^2)^c dv(z)
\\ &=& \sum_{k=1}^\infty\sum\limits_{|\alpha|=k}
k^{2s}|a_\alpha|^2\int_0^1 r^{2k}(1-r^2)^c r dr\int_{\partial
B_N}|\zeta^\alpha|^2 d\sigma(\zeta)
\\ &=& \sum_{k=1}^\infty\sum\limits_{|\alpha|=k}
k^{2s}|a_\alpha|^2
\frac{(N-1)!\alpha!}{(N-1+k)!}\cdot\frac{c!k!}{(c+k+1)!},
\end{eqnarray*}}which is equivalent to $$\sum_{k=1}^\infty\sum\limits_{|\alpha|=k}|a_\alpha|^2 \frac{(N-1)!\alpha!}{(N-1+k)!}\cdot k^{2s-c-1}$$
according to Stirling's formula. On the other hand, we  have
$$||f||^2_{H^2(\beta)}=\sum\limits_{k=0}^{\infty}||f_k||^2(k+1)^{2\nu}=
 \sum_{k=0}^\infty\sum\limits_{|\alpha|=k}|a_\alpha|^2 \frac{(N-1)!\alpha!}{(N-1+k)!}\cdot (k+1)^{2\nu}.$$
Note that $||f||_{A^2_{c, s}}=||R^sf||_{A^2_c}+|f(0)|$. If
$2s-c-1=2\nu$ and $s\ge \nu$, i.e. $c=2s-2\nu-1$ with $c\ge -1$,
then the norm of the weighted Hardy space $H^2(\beta, B_N)$ is
comparable with that of the  Sobolev space $A^2_{2s-2\nu-1, s}$. \ \
$\Box$

Thus,  using Lemma 4.1,  we have the following result for   weighted
Hardy spaces.
\\ \\ {\bf Proposition 4.2.}\begin{em}  Let $\varphi$ be a linear fractional
self-map of $B_N$. Then $C_\varphi$ is bounded on the weighted Hardy
space $H^2(\beta, B_N)$ with the weight $\beta(k)=(k+1)^v$, where
$\nu$ is a real number.\end{em}
\\ \par
Therefore, if we can calculate the spectral radii  of linear
fractional composition operators  on $H^2(\beta, B_N)$ with the
weight $\beta(k)=(k+1)^v$, one subject of this topic shall be  to
characterize  the spectra of these composition  operators.

\small DEPARTMENT OF MATHEMATICS, TONGJI UNIVERSITY,

SHANGHAI 200092, CHINA

DEPARTMENT OF APPLIED  MATHEMATICS,

SHANGHAI FINANCE UNIVERSITY,

SHANGHAI  201209, CHINA

{\it E-mail address:} liangying1231@163.com \\ \par DEPARTMENT OF
MATHEMATICS, TONGJI UNIVERSITY,

SHANGHAI 200092, CHINA

{\it E-mail address:}  zzzhhc@mail.tongji.edu.cn

\end{document}